\documentclass[12pt,a4paper]{article}
\usepackage{epsfig,latexsym,amsfonts,amssymb,amsmath,amscd,
epic,graphics,theorem}
\newtheorem{thm}{Theorem}
\newtheorem{lem}
{Lemma}
{Proposition}
{Claim}
\newtheorem{df}{Definition}
{Corollary}
\newtheorem{rem}
{Remark}
\newtheorem{que}
{Question}

{\theorembodyfont{\rmfamily}  \font\ncsc=cmcsc10
 \font\ntt=cmtt12

\begin{document}
\newcommand{\pperp}{\hbox{$\perp\hskip-6pt\perp$}}
\newcommand{\ssim}{\hbox{$\hskip-2pt\sim$}}
\newcommand{\aleq}{{\ \stackrel{3}{\le}\ }}
\newcommand{\ageq}{{\ \stackrel{3}{\ge}\ }}
\newcommand{\aeq}{{\ \stackrel{3}{=}\ }}
\newcommand{\bleq}{{\ \stackrel{n}{\le}\ }}
\newcommand{\bgeq}{{\ \stackrel{n}{\ge}\ }}
\newcommand{\beq}{{\ \stackrel{n}{=}\ }}
\newcommand{\cleq}{{\ \stackrel{2}{\le}\ }}
\newcommand{\cgeq}{{\ \stackrel{2}{\ge}\ }}
\newcommand{\ceq}{{\ \stackrel{2}{=}\ }}
\newcommand{\N}{{\mathbb N}}
\newcommand{\A}{{\mathbb A}}
\newcommand{\K}{{\mathbb K}}
\newcommand{\Z}{{\mathbb Z}}
\newcommand{\R}{{\mathbb R}}
\newcommand{\C}{{\mathbb C}}
\newcommand{\Q}{{\mathbb Q}}
\newcommand{\PP}{{\mathbb P}}
\newcommand{\Pic}{{\operatorname{Pic}}}\newcommand{\Sym}{{\operatorname{Sym}}}
\newcommand{\Aut}{{\operatorname{Aut}}}
\newcommand{\oeps}{{\overline\eps}}\newcommand{\idim}{{\operatorname{idim}}}
\newcommand{\oDel}{{\widetilde\Del}}
\newcommand{\cD}{{\mathcal D}}\newcommand{\mt}{{\operatorname{mt}}}
\newcommand{\ord}{{\operatorname{ord}}}\newcommand{\Id}{{\operatorname{Id}}}
\newcommand{\Span}{{\operatorname{Span}}}\newcommand{\As}{{\operatorname{As}}}
\newcommand{\Ker}{{\operatorname{Ker}}}\newcommand{\nd}{{\operatorname{nd}}}
\newcommand{\Ann}{{\operatorname{Ann}}}
\newcommand{\Fix}{{\operatorname{Fix}}}
\newcommand{\sign}{{\operatorname{sign}}}
\newcommand{\Tors}{{\operatorname{Tors}}}
\newcommand{\Ima}{{\operatorname{Im}}}
\newcommand{\oi}{{\overline i}}
\newcommand{\oj}{{\overline j}}
\newcommand{\ob}{{\overline b}}
\newcommand{\os}{{\overline s}}
\newcommand{\oa}{{\overline a}}
\newcommand{\oy}{{\overline y}}
\newcommand{\ow}{{\overline w}}
\newcommand{\ot}{{\overline t}}
\newcommand{\oz}{{\overline z}}
\newcommand{\eps}{{\varepsilon}}
\newcommand{\proofend}{\hfill$\Box$\bigskip}
\newcommand{\Int}{{\operatorname{Int}}}
\newcommand{\pr}{{\operatorname{pr}}}
\newcommand{\Hom}{{\operatorname{Hom}}}
\newcommand{\rk}{{\operatorname{rk}}}\newcommand{\Ev}{{\operatorname{Ev}}}
\newcommand{\im}{{\operatorname{Im}}}
\newcommand{\sk}{{\operatorname{sk}}}\newcommand{\DP}{{\operatorname{DP}}}
\newcommand{\const}{{\operatorname{const}}}
\newcommand{\Sing}{{\operatorname{Sing}}\hskip0.06cm}
\newcommand{\conj}{{\operatorname{Conj}}}
\newcommand{\Cl}{{\operatorname{Cl}}}
\newcommand{\Crit}{{\operatorname{Crit}}}
\newcommand{\Ch}{{\operatorname{Ch}}}
\newcommand{\discr}{{\operatorname{discr}}}
\newcommand{\Tor}{{\operatorname{Tor}}}
\newcommand{\Conj}{{\operatorname{Conj}}}
\newcommand{\vol}{{\operatorname{vol}}}
\newcommand{\defect}{{\operatorname{def}}}
\newcommand{\codim}{{\operatorname{codim}}}
\newcommand{\ov}{{\overline v}}
\newcommand{\ox}{{\overline{x}}}
\newcommand{\bw}{{\boldsymbol w}}
\newcommand{\bx}{{\boldsymbol x}}
\newcommand{\bd}{{\boldsymbol d}}
\newcommand{\bz}{{\boldsymbol z}}\newcommand{\bp}{{\boldsymbol p}}
\newcommand{\tet}{{\theta}}
\newcommand{\Del}{{\Delta}}
\newcommand{\bet}{{\beta}}
\newcommand{\kap}{{\varkappa}}
\newcommand{\del}{{\delta}}
\newcommand{\sig}{{\sigma}}
\newcommand{\alp}{{\alpha}}
\newcommand{\Sig}{{\Sigma}}
\newcommand{\Gam}{{\Gamma}}
\newcommand{\gam}{{\gamma}}
\newcommand{\Lam}{{\Lambda}}
\newcommand{\lam}{{\lambda}}
\newcommand{\SC}{{SC}}
\newcommand{\MC}{{MC}}
\newcommand{\nek}{{,...,}}
\newcommand{\cim}{{c_{\mbox{\rm im}}}}
\newcommand{\clM}{\tilde{M}}
\newcommand{\clV}{\bar{V}}

{\catcode`\@=11
\gdef\n@te#1#2{\leavevmode\vadjust{%
 {\setbox\z@\hbox to\z@{\strut#1}%
  \setbox\z@\hbox{\raise\dp\strutbox\box\z@}\ht\z@=\z@\dp\z@=\z@%
  #2\box\z@}}}
\gdef\leftnote#1{\n@te{\hss#1\quad}{}}
\gdef\rightnote#1{\n@te{\quad\kern-\leftskip#1\hss}{\moveright\hsize}}
\gdef\?{\FN@\qumark}
\gdef\qumark{\ifx\next"\DN@"##1"{\leftnote{\rm##1}}\else
 \DN@{\leftnote{\rm??}}\fi{\rm??}\next@}}
\def\mnote#1{\leftnote{\vbox{\hsize=2.5truecm\footnotesize
\noindent #1}}}

\title{On rigid plane curves}
\author{Victor Kulikov and Eugenii Shustin}
\date{}

\maketitle

\begin{abstract}
In the article, we exhibit a series of new examples of rigid plane curves, that is,
curves, whose collection of singularities determines them almost uniquely up to a projective transformation of the plane.
\end{abstract}

{\small

}

\section*{Introduction}
\subsection{Background and motivation} We work over the complex field, though most of results can be stated over any
algebraically closed field of characteristic zero.

The space ${\mathcal C}_d$ of plane curves of degree $d$ can be identified with $\mathbb P^{d(d+3)/2}$.
It has a natural {\bf equisingular stratification} with the strata determined by
the collection of degrees and multiplicities of irreducible components and by the collection of topological
singularity types of the considered curves (see \cite{Z1,Z2}; below, the strata will be called the {\it families
of equisingular curves}). Properties of this stratification have been studied by algebraic geometers since 19th
century, attracting attention of leading experts like Zeuthen, Severi, Segre, Zariski and others (see, for example,
\cite{GLS} for a modern survey in this area).

In this paper we focus on the minimal equisingular families,
that is, those which are formed by reduced curves and contain only finitely many orbits
of the action of the group of projective transformations of the plane $\Aut(\mathbb P^2)$.
The curves belonging to these families
are called {\bf rigid} curves; the corresponding families
we also call {\bf rigid} (see Definition \ref{rigid} below).
The study of rigid curves is motivated by their appearance in several important problems. First of all,
finite coverings of the projective plane branched along rigid curves are used to obtain
examples of rigid, so called, Miyaoka -- Yau surfaces (see, for example, \cite{Hir}, \cite{K-H}).
Next (cf. \cite{Pa}), celebrated Belyi's Theorem \cite{Be} says that each projective curve defined over
$\overline{\mathbb Q}$ can be represented as a finite covering of $\mathbb P^1$ branched at three points.
Note that any three points in $\mathbb P^1$ are rigid in the sense of Definition given below.
Therefore we can hope that for any field $F$ of transcendence degree two over $\overline{\mathbb Q}$
there are a projective model $X$ defined over $\overline{\mathbb Q}$
with the field of rational functions $\overline{\mathbb Q}(X)\simeq F$ and a finite morphism
$f:X\to \mathbb P^2$ branched along a rigid plane curve. There are interesting relations to the geometry
of line arrangements (see \cite{Hir}) and to rational cuspidal curves (see \cite{ZO}).

The goal of our note is to exhibit examples of rigid curves of any degree and any genus, and rigid families
covered by arbitrarily many orbits of the $\Aut(\mathbb P^2)$-action.

\subsection{Definitions and main results}
Throughout the paper we consider isolated plane curve singular points up to topological equivalence,
briefly calling any class of topologically equivalent singular points a {\it singularity type}.
Given a singularity type $S$, the number of irreducible components of singular curve germs of type $S$ is an invariant, which
we denote $m_S$.

Cardinality of a finite set $F$ will be denoted by $|F|$.

\subsubsection{Irreducible rigid curves}
Let $S_1,...,S_r$ be a sequence of distinct singularity types, $n_1,...,n_r$ a sequence of positive integers, $r\ge1$.
Introduce the formal sum
${\mathbf S}=\sum_{i=1}^rn_iS_i$.
Given a positive integer $d$, denote by $V(d;g;{\mathbf S})$ the
(equisingular) family of reduced, irreducible plane curves of degree $d$
having precisely $\sum_{i=1}^rn_i$ singular points and such that $n_i$ singular points are of type $S_i$, $i=1,...,r$.
Here $$g=\frac{(d-1)(d-2)}{2}-\sum_{i=1}^rn_i\delta(S_i)$$ is the geometric genus of the considered curves.

Such a family is a locally closed union of quasiprojective subvarieties of ${\mathcal C}_d$ (cf. \cite{GLS}).
It is, of course, invariant with respect to the action of $\Aut(\PP^2)$, and hence consists of entire orbits of
the $\Aut(\PP^2)$-action.

\begin{df} \label{rigid-irr} We say that a non-empty equisingular family $V(d,g,{\mathbf S})$ is {\bf $k$-rigid} if it is the
union of $k$ distinct orbits of the $\Aut(\PP^2)$-action in ${\mathcal C}_d$
for some $k\in \N$. If $k=1$ then we say that $V$ is {\bf strictly rigid.}
The curves belonging to a $k$-rigid family of plane curves are called {\bf rigid.}  \end{df}

\subsubsection{Reducible rigid curves}
Considering reduced, irreducible curves, we, first, introduce families of reducible curves
with numbered components, then identify families obtained from each other by permutation of components.

Let ${\bf d}=(d_1,\dots, d_N)$ and ${\bf g}=(g_1,\dots, g_N)$ be two collections of integers, $d_i\geq 1$ and
$g_i\geq 0$ for $i=1,\dots, N$. To encode the distribution of singularity types among
components and the distribution of local branches centered at singular points that are intersection points
of components, we do the following.
For a fixed singularity type $S$ and fixed $N$ denote by ${\bf J}_S=\{ J_{S,k}\}$ the set of all
nonempty subsets $J_{S,k}$ of $\{1,2,...,N\}$, $1\leq k\leq\sum_{j=1}^{m_S} \binom{N}{j}$, such that $|J_{S,k}|\leq m_{S}$.
Let $V({\bf d};{\bf g}; \sum_{\{ S_j\}}\sum_{\{ J_{S_j,k}\}}n_{J_{S_j,k}}S_j)$  be the family
of plane reduced curves $\overline C=C_1\cup\dots\cup C_N\subset \PP^2$ such that $C_i$ are irreducible curves
of degree $\deg C_i=d_i$ and genus $g_i$, and for each type $S_j$ of plane singularities the intersection
$\cap_{i\in J_{S_j,k}}C_i$ contains exactly $n_{J_{S_j,k}}$ singular points of $\overline C$ of the type $S_j$ which do
not lie in $C_l$ for $l\not\in J_{S_j,k}$. This is a locally closed union of quasiprojective subvarieties of ${\mathcal C}_d$ (cf. \cite{GLS}).

The sum ${\bf S}=\sum_{\{ S_j\}}\sum_{\{ J_{S_j,k}\}}n_{J_{S_j,k}}S_j$
is called the {\bf singularity type} of the curves $\overline C\in V({\bf d};{\bf g};{\bf S})$.
We identify the families $V({\bf d};{\bf g};{\bf S})$ obtained by permutations of
the curves $C_1,\dots,C_N$ and compatible permutations of ${\bf S}$.

A singularity type ${\bf S}$ splits into two parts, ${\bf S}={\bf S}^{ess}+{\bf S}^{non-ess}$, as follows.
For fixed ${\bf d}$ and ${\bf g}$, we say that ${\bf S}^{ess}$ is an {\bf essential part} of the singularity
type ${\bf S}$ (and resp. ${\bf S}^{non-ess}$ is a {\bf non-essential part} of the singularity type ${\bf S}$)
if the family  $V({\bf d;\bf g; \bf S})$ is determined uniquely by ${\bf d}$, ${\bf g}$, and the property that
the curves $\overline C$ have the singularities ${\bf S}^{ess}$ among the
all singularities of $\overline C$. If ${\bf S}^{ess}$ is an essential part of a singularity type ${\bf S}$,
then we will use notation $V({\bf d;\bf g;\bf S}^{ess}+\dots)$ to denote the family $V({\bf d;\bf g;\bf S})$.
\begin{df} \label{rigid} We say that a family $V=V({\bf d;\bf g; \bf S})$ is {\bf $k$-rigid} if it is the
union of $k$ distinct orbits of the $\Aut(\PP^2)$-action in ${\mathcal C}_d$
for some $k\in \N$. If $k=1$ then we say that $V$ is {\bf strictly rigid.}
The curves belonging to a $k$-rigid family of plane curves are called {\bf rigid.}  \end{df}

Note that the number of irreducible components of a $k$-rigid
family $V({\bf d;\bf g;\bf S})$ is less or equal $k$
(since it can (and does) happen that some orbits  can lie in the
closure of another) and, in particular, $V({\bf d;\bf g;\bf S})$ is
irreducible if it is strictly rigid. Note also that if a family $V({\bf d};{\bf g};{\bf S})$ is rigid, then
\begin{equation} \label{dim} \dim V({\bf d};{\bf g};{\bf S})\leq \dim PGL(\mathbb C,3)=8.\end{equation}

\subsubsection{Main results}
Our results are as follows:
\begin{itemize}
\item in Theorem \ref{d4}, the complete list of rigid curves of degree $\leq 4$ is given;
\item in Theorem \ref{cl2}, we give an infinite series of examples of strictly rigid families of
irreducible rational curves $V(d;0;{\bf S})$;
\item in Theorem \ref{t-oe1}, for each $g\geq 1$, we prove the existence of strictly rigid irreducible plane
curves of genus $g$;
\item examples of irreducible $2$-rigid families of irreducible curves are given in
Theorems \ref{add2} and \ref{add3}, and Theorem \ref{rigit} provides examples  of  $k$-rigid families
$V({\bf d};{\bf g};{\bf S})$ consisting of $k$ irreducible components for each $k\in \N$.
\end{itemize}

We do not know answers to the following questions, which seem to be interesting:

\begin{que}
Do there exist  irreducible $k$-rigid
families $V({\bf d};{\bf g};{\bf S})$ with $k> 2${\rm ?}

Do there exist irreducible $2$-rigid families $V(d;g;{\bf S})$ with $g\geq 1${\rm ?}\end{que}

Through the paper, we use the following notations for
singularity types of plane curves:
\begin{itemize}
\item
$T_{m,n}$, $2\leq m\leq n$, is the type of singularity given by equation $x^m+y^{n}=0$;
but if $m=2$ then a singularity of type $T_{2,n}$, as usual, will be denoted by $A_{n-1}$ and the
singularities of types $T_{m,m}$ will be called {\it simple.}
\item
$T^m_{m,n}$, $2\leq m< n$, is the type of singularity
given by equation $y(x^m+y^{n})=0$.
\item
$T^n_{m,n}$, $2\leq m< n$, is the type of singularity
given by equation $x(x^m+y^{n})=0$.
\item
$T^{m,n}_{m,n}$, $1\leq m< n$, is the type of singularity
given by equation $xy(x^m+y^{n})=0$.
\end{itemize}

{\small
{\bf Acknowledgements}.
The first author has been supported
by grants of  RFBR 14-01-00160 and 15-01-02158, and by the Government of the Russian Federation within the
framework of the implementation of the 5-100 Programme Roadmap of the National Research University  Higher
School of Economics, AG Laboratory. The second author has been supported by
the Hermann-Minkowski-Minerva Center for Geometry at the Tel Aviv
University and by the German-Israeli
Foundation grant no. 1174-197.6/2011. Main ideas behind this work have appeared during
the visit of the second author to the Higher School of Economics (Moscow).
We are grateful to
these institutions for support and excellent working conditions.}
\section{Rigid  curves of small degree}
In the following Theorem, we provide the complete list of rigid reduced curves of degree $\leq 4$.

\begin{thm} \label{d4} Let $\overline C$ be a rigid reduced curve of degree $\leq 4$. Then $\overline C$ belongs to
one of the following families:
\begin{itemize}
\item[$(I)$] strongly rigid families:
\begin{itemize}
\item[$(I_1)$] $V(1;0;\emptyset)$;
\item[$(I_2)$] $V((1,1);(0,0);A_1)$;
\item[$(I_3)$] $V(2;0;\emptyset)$;
\item[$(I_4)$] $V((1,1,1);(0,0,0);3A_1)$;
\item[$(I_5)$] $V((1,1,1);(0,0,0);T_{3,3})$;
\item[$(I_6)$] $V((2,1);(0,0);2A_1)$;
\item[$(I_7)$] $V((2,1);(0,0);A_3)$;
\item[$(I_8)$] $V(3;0;A_1)$;
\item[$(I_9)$] $V(3;0;A_2)$;
\item[$(I_{10})$] $V((1,1,1,1);(0,0,0,0);6A_1)$;
\item[$(I_{11})$] $V((1,1,1,1);(0,0,0,0);T_{3,3}+3A_1)$;
\item[$(I_{12})$] $V((2,1,1);(0,0,0);2A_3+A_1)$;
\item[$(I_{13})$] $V((2,1,1);(0,0,0);A_3+3A_1)$;
\item[$(I_{14})$] $V((2,1,1);(0,0,0);T_{2,4}^2+A_1)$;
\item[$(I_{15})$] $V((2,2);(0,0);A_5+A_1)$;
\item[$(I_{16})$] $V((2,2);(0,0);A_7)$;
\item[$(I_{17})$] $V((3,1);(0,0);A_5+A_1)$;
\item[$(I_{18})$] $V((3,1);(0,0);T_{2,4}^2)$;
\item[$(I_{19})$] $V((3,1);(0,0);T_{2,3}^3)$;
\item[$(I_{20})$] $V((3,1);(0,0);A_{2}+A_5)$;
\item[$(I_{21})$] $V((3,1);(0,0);A_3+A_2+A_1)$;
\item[$(I_{22})$] $V(4;0;3A_2)$;
\item[$(I_{23})$] $V(4;0;A_4+A_2)$;
\item[$(I_{24})$] $V(4;0;A_{6})$;\end{itemize}
\item[$(II)$] irreducible $2$-rigid families:
\begin{itemize}\item[$(II_{1})$] $V((3,1);(0,0);T_{2,3}^2+A_1)$;
\item[$(II_{2})$] $V(4;0;T_{3,4})$.
\end{itemize}
\end{itemize}
\end{thm}
{\bf Proof.}
It is well known that if $\deg \overline C\leq 3$ then only smooth cubics are not rigid. All other families $V({\bf d};{\bf g};
{\bf S})$ of curves of degree $\leq 3$ are listed in $(I_1)$ -- $(I_9)$. Therefore we assume below that $\deg \overline C=4$.

Again, it is well known that if $\overline C$ consists of four lines or two lines and a quadric, then $\overline C$ is not rigid
if and only if $\overline C$ consists of four lines having a common point or $\overline C$ consists of a quadric $Q$ and two lines
in general position with respect to $Q$. The rigid families in the case when $\overline C$ consists of four lines or two lines and
a quadric  is listed in $(I_{10})$ -- $(I_{14})$.

Consider the case when $\overline C$ consists of two irreducible components: either $\overline C=Q_0\cup Q_1$, where $Q_0$ and
$Q_1$ are smooth quadrics, or  $\overline C=C\cup L$, where $C$ is a cubic and $L$ is a line.

Consider the case when $\overline C_1=Q_0\cup Q_1\in V((2,2);(0,0);{\bf S})$. We have ${\bf S}=m_1A_1+m_3A_3+m_5A_5+m_7A_7$,
where $m_1+2m_3+3m_5+4m_7=4$.  The singularity type ${\bf S}$ consists of $k=m_1+m_3+m_5+m_7$, $1\leq k\leq 4$, singular points
and if  $m_1=4$, that is, ${\bf S}=4A_1$, then $V((2,2);(0,0);4A_1)$ is not rigid by inequality \eqref{dim}, since
$\dim V((2,2);(0,0);4A_1)=10$. Similarly, if $m_1=2$, that is, ${\bf S}=A_3+2A_1$, then $V((2,2);(0,0);A_3+2A_1)$ is
not rigid, since $\dim V((2,2);(0,0);A_3+2A_1)=9$. So, we can assume that $m_1\leq 1$.

Consider the pencil of quadrics $Q_{\lambda}$ defined by quadrics $Q_0$ and $Q_1$. It is easy to see that in the cases
$m_1=m_5=1$ or $m_7=1$  it contains the unique degenerate element $Q_{\infty}$ consisting of two (coinciding if $m_7=1$)
lines $L_1\cup L_2$ and in the case $m_3=2$ it contains two degenerate elements one of which, $Q_{\infty}$, consists of
two coinciding lines $L_1=L_2$ and the other one consists of two different lines.
Let $f(z_1,z_2,z_3)=0$ be an equation of the quadric $Q_0$ and $l_i(z_1,z_2,z_3)=0$, $i=1,2$, be an equation of the line
$L_i$. Then, without less of generality, we can assume that $$F_{\lambda}(z_1,z_2,z_3):= f(z_1,z_2,z_3)+\lambda
l_1(z_1,z_2,z_3)l_2(z_1,z_2,z_3)=0$$ is the equation of $Q_{\lambda}$. Moreover, applying a
projective transformation, we can assume that $f(z_1,z_2,z_3)=z_1^2-z_2z_3$ and $l_1(z_1,z_2,z_3)=z_3$, $l_1(z_1,z_2,z_3)=
z_1$ in the case $m_1=m_5=1$;  $l_1(z_1,z_2,z_3)=l_2(z_1,z_2,z_3)=z_1$ in the case $m_3=2$; and $l_1(z_1,z_2,z_3)=l_2(z_1,
z_2,z_3)=z_3$ in the case $m_7=1$. Therefore we can assume that
\begin{equation}\label{pen1} F_{\lambda}(z_1,z_2,z_3)= (z_1^2-z_2z_3)+\lambda z_1z_3
\end{equation}
in the case $m_1=m_5=1$,
\begin{equation}\label{pen3} F_{\lambda}(z_1,z_2,z_3)= (z_1^2-z_2z_3)+\lambda z_3^2
\end{equation}
in the case $m_7=1$, and
\begin{equation}\label{pen2} F_{\lambda}(z_1,z_2,z_3)= (z_1^2-z_2z_3)+\lambda z_1^2
\end{equation}
in the case $m_3=2$.

Note that if $m_3=2$, then the degenerate element $Q_{-1}$ of the pencil $Q_{\lambda}$ consists of two lines $z_2=0$ and $z_3=0$
(see equation \eqref{pen2}).

In the case $m_1=m_5=1$ (case $(I_{15}))$ the strong rigidity of the curve $\overline C_1$ follows from equality $\overline
C_{\lambda}=Q_0\cup Q_{\lambda}=h_{\lambda}(\overline C_1)$ for each $\lambda\in \mathbb C^*$, where the automorphism
$h_{\lambda}$ acts as follows: $h_{\lambda}(z_1:z_2:z_3)=(\lambda z_1: \lambda^2 z_2:z_3)$; and in the case $m_7=1$
(case $(I_{16}))$ the strong rigidity of the curve $\overline C_1$ follows from equality $\overline C_{\lambda^2}=Q_0
\cup Q_{\lambda^2}=h_{\lambda}(\overline C_1)$ for each $\lambda\in \mathbb C^*$.

Let us show that $\overline C_1$ is not rigid if $m_3=2$. Indeed, assume that $\overline C_1$ is rigid.
Then for any two elements  $\lambda_{1},\lambda_{2}\in \mathbb C^*\setminus \{ -1\})$,  $\lambda_{1}\neq \lambda_{2}$, there is a
projective transformation $h_{{\lambda}_1,{\lambda}_2}\in\Aut(\mathbb P^2)$ such that $h_{{\lambda}_1,
{\lambda}_2}(\overline C_{{\lambda}_1})=\overline C_{{\lambda}_2}$, where
$\overline C_{\lambda }=Q_{0}\cup Q_{\lambda}$. The automorphism  $h_{{\lambda}_1,{\lambda}_2}$ leaves
invariant the pencil $Q_{\lambda}$. In particular, it leaves invariant the singular elements $Q_{\infty}$
and $Q_{-1}$ of the pencil $Q_{\lambda}$, $$h_{{\lambda}_1,{\lambda}_2}(Q_{\infty})=Q_{\infty}\, \, \,
\text{and}\, \, \, h_{{\lambda}_1,{\lambda}_2}(Q_{-1})=Q_{-1}.$$
Let us show that $h_{\lambda_1,\lambda_2}(Q_0)\neq Q_0$ for $\lambda_1\neq \lambda_2$. Indeed, if
$h_{\lambda_1,\lambda_2}(Q_0)= Q_0$, then in the case $m_3=2$ we have either $h^*_{\lambda_1,\lambda_2}(z_1)
= az_1$, $h^*_{\lambda_1,\lambda_2}(z_2)= bz_2$, and $h^*_{\lambda_1,\lambda_2}(z_3)= cz_3$ or $h^*_{\lambda_1,
\lambda_2}(z_1)= az_1$, $h^*_{\lambda_1,\lambda_2}(z_2)= bz_3$, and $h^*_{\lambda_1,\lambda_2}(z_3)= cz_2$ for
some $a,b,c$ such that $a^2=bc$, since $h_{\lambda_1,\lambda_2}(Q_{\infty})= Q_{\infty}$ and $h_{\lambda_1,\lambda_2}
(Q_{-1})= Q_{-1}$. Therefore $h_{\lambda_1,\lambda_2}(Q_{\lambda})= Q_{\lambda}$ for all $\lambda \in \mathbb C$
that is possible only if $\lambda_1=\lambda_2$. Therefore we must have $h_{\lambda_1,\lambda_2}(Q_{0})= Q_{\lambda_2}$
and $h_{\lambda_1,\lambda_2}(Q_{\lambda_1})= Q_{0}$. But, for three pairwise different $\lambda_1,\lambda_2,\lambda_
3\in \mathbb C^*\setminus \{ -1\}$ we obtain that
$h_{\lambda_1,\lambda_3}(Q_0)=h_{\lambda_2,\lambda_3}\circ h_{\lambda_1,\lambda_2}(Q_0)=Q_0$ and hence $\lambda_1=\lambda_3$. Contradiction.

Consider the case when $\overline C=C\cup L$, where $C$ is a cubic and $L$ is a line. It is easy to see that
$\overline C$ is rigid only if $C$ is a rational curve. Therefore we have two cases: $C$ is a nodal cubic or $C$ is a cuspidal cubic.

Let $C$ be a nodal cubic. Then (the case $(I_8)$) $V(3;0;A_1)$ is strongly rigid and $\dim V(3;0;A_1)=8$.
Therefore $L$ must be a "very special line" with respect to $C$, that is, $L$ is either the tangent line of one of
 two branches of the node of $C$ or $L$ is the tangent line of $C$ at a flex point of $C$. In both two cases it is easy to see that $\overline C$ is strictly rigid and we have two cases: $(I_{17})$ -- $(I_{18})$.

If $C$ is a cuspidal cubic, then $L$ is a tangent line to $C$ or it passes through the cusp of $C$, since $V(3;0;A_2)$ is strongly rigid and $\dim V(3;0;A_2)=7$. If $L$ is a tangent line to $C$ at its cusp, then we have case ($I_{19}$) and if $L$ is the tangent line at the flex point of $C$, then we have case $(I_{20})$. It is easy to see that in both cases $\overline C$ is strongly rigid.

If $L$ is a tangent line to $C$ and it does not pass through the cusp of $C$, let us choose homogeneous coordinates $(z_1:z_2:z_3)$ such that $z_3=0$ is an equation of $L$, $z_1=0$ is an equation of the line tangent to $C$ at its cusp, and the line given by equation $z_2=0$ passes through the cusp of $C$ and the tangent point of $L$ and $C$. Then there is a parametrization of $C$ in $\mathbb P^2$ of the following form
\begin{equation} \label{parx} z_1=t^3,\quad z_2=t^2,\quad z_3=1+at,\end{equation}
where $a\neq 0$. Denote by $C_a$ a curve in $\mathbb P^2$ given by parametrization \eqref{parx}. Now, the strong rigidity of $V((3,1);(0,0);A_3+A_2+A_1)$  (case $(I_{21})$) follows from equality $h(C_a\cup L)=C_1\cup L$, where the linear transformation $h\in \Aut(\mathbb P^2)$ acts as follows, $h(z_1:z_2:z_3)=(z_1:a^{-1}z_2:a^{-3}z_3)$.

Let $L$ be not the tangent line at the cusp of $C$ and it pass through the cusp of $C$. Then if $z_2=0$ is an equation of $L$, $z_1=0$ is the line tangent to $C$ at its cusp, and the line given by equation $z_3=0$ is tangent to $C$ at its nonsingular  point $p$ of the intersection $C\cap L$, then again we can assume that $C$ is given by parametrization 
\eqref{parx} and there are two possibilities: either $a=0$ or $a\neq 0$. Denote by $C_a$ a curve given by parametrization \eqref{parx}. It is easy to see that $p\in C_0\cap L$ is the flex point of $C_0$ if $a=0$ and $p\in C_a\cap L$ is not the flex point of $C_a$ if $a\neq 0$. Therefore there is not a linear transformation $h\in \Aut(\mathbb P^2)$ such that $h(C_0\cup L)=C_a\cup L$. On the other hand, the linear transformation $h:(z_1:z_2:z_3)\rightsquigarrow (z_1:a^{-1}z_2:a^{-3}z_3)$ sends $C_a\cup L$ to $C_1\cup L$. Therefore $V((3,1);(0,0), T_{2,3}^2)$ (case $(II_1)$) is an irreducible $2$-rigid family.

Now, consider the case when $\overline C$ is an irreducible curve of degree four. All possible type of singularities ${\bf S}$ of irreducible curves of degree four is given in \cite{Wal}:
$$ {\bf S}=m_1A_1+m_2A_2+\dots m_6A_6+m_7T_{3,3}+m_8T_{2,3}^2+m_9T_{3,4},$$
where $m_i$, $1\leq i\leq 9$, are non-negative integers satisfying the following inequality:
\begin{equation} \label{type}
m:=m_1+m_2+2m_3+2m_4+3m_5+3m_6+3m_7+3m_8+3m_9\leq 3.\end{equation}
By \cite[Theorem 2]{Sh}, we have
\begin{equation} \label{dimV} \begin{array}{l} \dim V(4;3-m;m_1A_1+m_2A_2+\dots+ m_6A_6+m_7T_{3,3}+m_8T_{2,3}^2+m_9T_{3,4})= \\
14 - (m_1+2m_2+3m_3+4m_4+5m_5+6m_6+4m_7+5m_8+6m_9).
\end{array}\end{equation}
Inequality \eqref{dim} applied to \eqref{dimV} gives rise to the inequality
\begin{equation} \label{di} 14 - (m_1+2m_2+3m_3+4m_4+5m_5+6m_6+4m_7+5m_8+6m_9)\leq 8. \end{equation}
It is easy to see that inequalities \eqref{type} and \eqref{di} have only the following non-negative integer solutions:
\begin{itemize} \item[$1)$] $m_2=3$ and $m_i=0$ for $i\neq 2$,
\item[$2)$] $m_2=m_4=1$ and $m_i=0$ for $i\neq 2,4$,
\item[$3)$] $m_6=1$ and $m_i=0$ for $i\neq 6$,
\item[$4)$] $m_9=1$ and $m_i=0$ for $i\leq 8$. \end{itemize}
Therefore only the following irreducible curves $\overline C$ of degree four can be rigid:
either $\overline C\in V(4;0;3A_2)$, or $\overline C\in V(4;0;A_2+A_4)$, or $\overline C\in V(4;0;A_6)$, or $\overline C\in V(4;0;T_{3,4})$.

The strict rigidity of a curve $\overline C\in V(4;0;3A_2)$ (case $(I_{22})$) is well known.

To show that a curve $\overline C\in V(4;0;A_2+A_4)$  (case $(I_{23})$) is strictly rigid, let us choose
a non-homogeneous coordinate $t$ in $\PP^1$ and homogeneous coordinates $(z_1:z_2:z_3)$ in $\PP^2$ such that $f(t=0)=(0:0:1)$ is the singular point of $\overline C$ of type $A_4$ and the line $z_1=0$ is tangent  to $\overline C$ at this point,  the image $f(t=\infty)=(1:0:0)$ is the singular point of $\overline C$, and the line $z_3=0$ is tangent to $C$ at this point. Then (after a suitable linear change of the coordinate $t$) the morphism $f$ is given by
\begin{equation} \label{nu1} z_1=t^4, \, \,
z_2=t^2, \, \, z_3=t+1, \end{equation}
that is, $\overline C$ is strictly rigid.

To show that  a curve $\overline C\in V(4;0;A_6)$ (case $(I_{24})$) is strictly rigid, note that, by Pl\"ucker's formulas, $\overline C$ must have a flex point. As above, let $f:\PP^1\to\PP^2$
be a morphism such that $\overline C=f(\PP^1)$. We can choose a non-homogeneous coordinate $t$ in $\PP^1$ and homogeneous coordinates $(z_1:z_2:z_3)$ in $\PP^2$ such that $f(t=0)=(0:0:1)$ is the singular point of $C$,
the line $z_1=0$ is the tangent line to $C$ at its singular point,  the image $f(t=\infty)=(1:0:0)$ is a flex point of $C$, and, moreover, the line $z_3=0$ is the tangent line to $C$ at
this point. Then (after a suitable linear change of the coordinate $t$) the morphism $f$ is given by
\begin{equation} \label{nu2} z_1=t^4, \, \,
z_2=t^2(t-1), \, \, z_3=F(t):=a_1t+a_2, \end{equation}
where $a_2\neq 0$. It is easy to see that the polynomial $F(t)-a_2(t-1)^2$ must be divisible
by $t^2$, since $C$ has a singularity of type $A_6$ at $(0:0:1)$.
Therefore we have $a_1=-2a_2$ and we can put $a_2=1$, that is, $V(4;0;A_6)$ is strictly rigid.

Let us show that the family $V(4;0;T_{3,4})$ (case $(II_{2})$) is irreducible and $2$-rigid.
Indeed, let $\overline C\in V(4;0;T_{3,4})$, then, by Pl\"ucker's formulas, $\overline C$ must have a flex point. Let $f:\PP^1\to\PP^2$ be a morphism such that $\overline C=f(\PP^1)$. We can choose
non-homogeneous coordinate $t$ in $\PP^1$ and homogeneous coordinates $(z_1:z_2:z_3)$
in $\PP^2$ such that $f(t=0)=(0:0:1)$ is the singular point of $\overline C$, the line $z_1=0$ is the
tangent line to $\overline C$ at its singular point,  the image $f(t=\infty)=(1:0:0)$ is a flex point of $\overline C$, and, moreover, the line $z_3=0$ is the tangent line to $\overline C$ at this point. Then (maybe, after
change of the coordinate $t$) the morphism $f$ is given by
\begin{equation} \label{nu3} z_1=t^4, \, \,
z_2=t^3, \, \, z_3=at+1. \end{equation}
Therefore $V(4;0;T_{3,4})$ is an irreducible family.

As in case $(II_1)$, there are two possibilities: either $a=0$ or $a\neq 0$.

Denote by $\overline C_a\in V(4;0;T_{3,4})$ a curve given by equation \eqref{nu3}. It is easy to see that
$\overline C_0$ has the unique flex point, namely, $p_1=(1:0:0)$, and one can check that
the curve $\overline C_1$ has two flex points, $p_1$ and $p_2=f(t=-2)=(16:-8:-1)$. Hence there is no
projective transformation $h\in \Aut(\PP^2)$ such that $\overline C_0=h(\overline C_1)$.

If we perform the change $t_1=at$, where $a\neq 0$, then we get $h(\overline C_1)=\overline C_a$, where
$$h(z_1:z_2:z_3)=(a^{4}z_1:a^{3}z_2:z_3)\ ,$$
that is, $V(4;0;T_{3,4})$ is $2$-rigid. \proofend

\begin{rem} {\rm According to Theorem \ref{d4}, the families
$V((1,1,1,1);(0,0,0,0);T_{4,4})$ and $V((2,2);(0,0); 2A_3)$ are not rigid and it is easy to see that $$\dim V((1,1,1,1);(0,0,0,0);T_{4,4})=6\, \, \, \,  \text{and} \, \, \, \, \dim V((2,2);(0,0); 2A_3)=8,$$ that is, the
inequality \eqref{dim} is not sufficient for $V({\bf d};{\bf g};{\bf S})$ to be a rigid family.} \end{rem}

\section{Strictly rigid rational curves of degree $\geq 5$}
\begin{thm} \label{cl2} For each $n\geq 2$ the family $V_n=V(2n; 0;3T_{n,n+1}+\dots)$ is strictly rigid. The non-essential part of singularities of $C\in V_n$ consists of simple singularities.
\end{thm}

To prove Theorem \ref{cl2} we need in the following
\begin{thm} \label{cl3} Let ${\bf d}_n=(n,1,1,1)$, ${\bf g}=(0,0,0,0)$, and $n_{(1,2)}=n_{(1,3)}=n_{(1,4)}=1$. Then the family
$$\hat{V}_n=V({\bf d}_n;{\bf g};(n_{(1,2)}+n_{(1,3)}+ n_{(1,4)})A_{2n-1}+\dots)$$
is strictly rigid for each $n\geq 2$. The non-essential part of singularities of $\overline C\in \hat{V}_n$ consists of simple singularities.
\end{thm}

We prove Theorems \ref{cl2} and \ref{cl3} simultaneously. \newline {\bf Proof.} Theorem \ref{cl2} in the case $n=2$ and Theorem \ref{cl3} in the case $n\leq 3$ are well known, and we prove the general case by induction on $n$.

Assume that for $n\leq n_0$ Theorem \ref{cl2} is true. Then the number of virtual cusps of $C\in V_n$ ($n
\leq n_0$) is equal to $3(n-1)$ and the number of its virtual nodes is
equal to $2(n-1)(n-2)$ (see the definition of the numbers of virtual cusps and virtual nodes, for example, in \cite{Ku}). Let $\hat{C}$ be the dual curve to $C$. By Pl\"ucker's formulas, we have $\deg \hat{C}=n+1$ and the number of virtual cusps of the curve $\hat{C}$ is zero.
Therefore the irreducible branches of the singular points of
$\hat{C}$ are smooth and hence $\hat C$ has only simple singularities, since $C$ has the only simple
singularities and three singularities of type $T_{n,n+1}$.  Note also that $\hat{C}$ has three flex points, say $p_1$, $p_2$, and $p_3$.
Let $L_i$, $i=1,2,3$, be the tangent line to $\hat{C}$ at the point $p_i$. It is easy to see that
$$\overline C=\hat{C}\cup L_1\cup L_2\cup L_3\in \hat{V}_{n+1}=V({\bf d}_{n+1};
{\bf g};(n_{(1,2)}+n_{(1,3)}+n_{(1,4)})A_{2n+1}+\dots).$$
Moreover, if $n\leq n_0$ and the curve $\overline C_1=C_1\cup L'_1\cup L'_2\cup L'_3$
belongs to $\hat{V}_{n+1}$, then the curve $C_1$ must have only simple singularities, since the dual curve $\hat{C}_1$ belongs to $V_{n}$.
Therefore the assumption that Theorem \ref{cl2} is true for $n\leq n_0$ implies
the statement of Theorem \ref{cl3} for $n\leq n_0+1$.\proofend

\begin{rem} \label{rem} {\rm For $n\leq n_0+1$, it follows from strong rigidity that if $\overline C_i=C_i\cup L_1\cup L_2\cup L_3\in \hat{V}_{n}$, $i=1,2$, are two curves such that $C_1$ and $C_2$ have two common flex points, say $p_1$ and $p_2$, then their third flex points $p_3=C_1\cap L_3$ and $q_3=C_2\cap L_3$
should coincide,  $p_3=q_3$, and consequently $C_1=C_2$.} \end{rem}

To complete the proof of Theorems \ref{cl2} and \ref{cl3}, note that if $\overline C=C \cup L_1\cup L_2\cup L_3$
belongs to $\hat{V}_{n}$ for $n\leq n_0+1$, then $\sigma(C)\in V_{n+1}$, where $\sigma$ is a quadratic transformation of the plane with centers at the vertices of the triangle $L_1\cup L_2\cup L_3$. Therefore Theorem \ref{cl2} in the case $n=n_0+1$ follows from Theorem \ref{cl3} and Remark \ref{rem}. \proofend

The following theorem provides an infinite series of strictly rigid families parameterizing the unions of two rational curves.

\begin{thm} \label{add1}
The families $V((2n+1,1);(0,0);n_{(1,2)}(T_{2n-1,2n}^{2n}+A_1)+n_{(1)}A_{4n-2})$, where $n_{(1,2)}=n_{(1)}=1$, are strictly rigid for $n\geq 2$.
\end{thm}
{\bf Proof.} Let $$C\cup L\in V((2n+1,1);(0,0);(0,0);n_{(1,2)}(T_{2n-1,2n}^{2n}+A_1)+n_{(1)}A_{4n-2}).$$ Then
$C\in V(2n+1;0;T_{2n-1,2n}+A_{4n-2})$ and $L$ is the line tangent to $C$ at its singular point $p_1$ of type $T_{2n-1,2n}$.
Let $f:\PP^1\to\PP^2$ be a morphism such that $C=f(\PP^1)$. We can choose non-homogeneous coordinate $t$ in $\PP^1$ and
homogeneous coordinates $(z_1:z_2:z_3)$ in $\PP^2$ such that $f(t=\infty)=p_1=(1:0:0)$ is the singular point of $C$ of type
$T_{2n-1,2n}$ and the line $L$ is given by equation $z_3=0$,  $f(t=0)=p_2=(0:0:1)$ is the singular point of $C$ of type $A_{4n-2}$ and
the line $z_1=0$ is the tangent line to $C$ at $p_2$. Then (maybe, after change of
the coordinates in $\PP^2$ and $\PP^1$) the morphism $f$ is given by
\begin{equation} \label{par2} z_1=t^4P_{2n-3}(t), \, \,
z_2=t^{2}, \, \,z_3=t-1 ,\end{equation}
where $P_{2n-3}(t)=t^{2n-3}+\sum_{i=0}^{2n-4}a_it^{i}$.

To prove Theorem \ref{add1}, it suffices to show that the polynomial $P_{2n-3}(t)$ is defined uniquely by
the property that the point $p_2$ is a singular point of $C$ of type $A_{2(2n-1)}$.
For this put $x=z_2/z_3$, $y=z_1/z_3$ and $x_2=x$, $y_2=y/x^2$. Then the germ of
singularity $(C,0)\subset(\C^2,o)$, given by parametrization (at $t=0$):
$$x_2=t^{2}/(t-1), \, \, y_2=(t-1)P_{2n-3}(t),$$
has singularity  type $A_{2(2n-3)}$ at the point $o=(0,-a_0)$.

Consider a sequence of polynomials $A_k(t)=t^k+\sum_{i=0}^{k-1}a_{i,k}t^i$, $k\in \N$, where $A_1(t) = t+1$ and
\begin{equation} \label{seq}
\displaystyle A_{k+1} = \frac{t^2A_k(t)-A_k(1)}{t-1}  \end{equation}
for each $k\geq 1$.

Now, Theorem \ref{add1} follows from
\begin{lem} \label{l} Let a germ of singularity $(C,0)\subset(\C^2,o)$ of type $A_{2k}$, $k\geq 1$, be given by parametrization of the form (at $t=0$):
$$
x=t^{2}/(t-1), \, \, y=(t-1)P_{k}(t),$$
where $P_{k}(t)=t^k+\sum_{i=0}^{k-1}a_it^i$  and $o=(0,-a_0)$. Then the polynomial $P_{k}(t)$ is defined uniquely and, moreover, $P_k(t)=A_k(t)$.
\end{lem}
{\bf Proof.} If $k=1$ then it is easy to see that $P_1(t)=t+1$.
Assume that for $k<k_0$ Lemma is true and prove it in the case $k=k_0$.

If $(C,0)\subset(\C^2,o)$ is a germ of singularity of type $A_{2k_0}$, then the polynomial $(t-1)P_{k_0}(t)+a_0$  is divisible by $t^2$. Therefore  $(t-1)P_{k_0}(t)+a_0 =t^2Q_{k_0}(t)$, where $Q_{k_0}(t)=t^{k_0-1}+\sum_{i=0}^{k_0-2}b_it^i$ is a polynomial of degree $k_0-1\geq 1$ and the singularity given by
$$
x_1=t^{2}/(t-1), \, \, y_1=(t-1)Q_{k_0-1}(t)$$
is of type $A_{2(k_0-1)}$. Then, by assumption, $Q_{k_0-1}(t)=A_{k_0-1}(t)$ and hence
$$(t-1)P_{k_0}(t)+a_0 =t^2A_{k_0-1}(t).$$ In particular, $a_0=A_{k_0-1}(1)$, that is,
$$\displaystyle P_{k_0}(t) = \frac{t^2A_{k_0-1}(t)-A_{k_0-1}(1)}{t-1}.$$
\proofend

\section{Strictly rigid curves of positive genera}
The following theorem states that the strongly rigid family $V(4;0;3A_2)$
(see case $(I_{22})$ in Theorem \ref{d4}) is the first member in an infinite sequence of strictly rigid families.

\begin{thm}\label{t-oe1}
For any $n\ge2$, the family $V(2n,\left[\frac{n}{2}\right]-1,A_n+2T_{n,2n-1})$ is non-empty and
strictly rigid.
\end{thm}
{\bf Proof.}
Let $C$ be a plane curve of degree $2n$ with singularities $A_n+2T_{n,2n-1}$. Choose projective coordinates
$(z_1:z_2:z_3)$ so that the singular points of type $T_{n,2n-1}$ are located at
$(0:1:0)$ and $(0:0:1)$ with tangent lines $z_3=0$ and $z_2=0$, respectively, and the remaining
intersection points of $C$ with these lines are $(1:-1:0)$ and $(1:0:-1)$, respectively. Then,
in the affine coordinates $x=z_2/z_1$, $y=z_3/z_1$, the curve $C$ is given by an equation
\begin{equation}1+x+y+\sum_{i=1}^na_ix^iy^i=0,\quad a_n\ne0\ ,
\label{e-oe1}\end{equation}
the following lemma completes the proof of Theorem:

\begin{lem}\label{l-oe1}
There exists a unique curve in $\C^2$, given by equation \eqref{e-oe1} and having a singularity
of type $A_n$ in $(\C^*)^2$. It is irreducible and has genus $\left[\frac{n}{2}\right]-1$.
\end{lem}
{\bf Proof.} Note that a singularity of type $T_{n,2n-1}$ is analytically irreducible,
$\delta(T_{n,2n-1})=(n-1)^2$, and if a plane curve $C'$ has a singular point
of type $T_{n,2n-1}$ then $\deg C'\geq 2n-1$. Therefore the curve $C$ in question is
irreducible, since if a curve $C'$ of degree $2n-1$ has two singular points of type $T_{n,2n-1}$, then
$$g(C')\leq \frac{(2n-2)(2n-3)}{2}-2\delta(T_{n,2n-1})=-(n-1)<0,$$
a contradiction.

Under assumption that $C$ has no other singularities in
$\C^2$, the genus value follows from the fact that
$\delta(T_{n,2n-1})=(n-1)^2$ and $\delta(A_n)=\left[\frac{n+1}{2}\right]$.

To find a curve given by \eqref{e-oe1} and having a singularity of type $A_n$, substitute
$(x,y)$ for $(x/y,y)$ in \eqref{e-oe1} and multiply by $y$:
$$y^2+y\left(1+\sum_{i=1}^na_ix^i\right)+x=0\ .$$
This is a quadratic equation in $y$ with the discriminant
$$ \Delta(x)=\left( 1+\sum_{i=1}^na_ix^i\right)^2-4x\ .$$
The existence of a singularity of type $A_n$ on the considered curve in
$(\C^*)^2$ is equivalent to the condition that $\Delta(x)$ has a root of multiplicity $n+1$. We shall show that
this condition has a unique
solution $(a_1,...,a_n)$ such that $a_n\ne0$, and that for this solution
the other roots of $\Delta$ are simple. So, into the relation
\begin{equation}\Delta(x)=\left(1+\sum_{i=1}^na_ix^i\right)^2-4x\quad
\text{is divisible by}\quad (1+tx)^{n+1},\ t\ne0\ ,
\label{e-oe2}\end{equation}
we substitute $x=z^2$, $t=-\tau^2$, and then get an equivalent condition
\begin{equation}\left(1+\sum_{i=1}^na_iz^{2i}-2z\right)\left(1+\sum_{i=1}^na_iz^{2i}+2z\right)
\ \text{is divisible by}\ (1-\tau^2 z^2)^{n+1},\ \tau\ne0\ .\label{e-oe5}\end{equation}
Since the difference $4z$ of the factors in the former product in \eqref{e-oe5} is not
divisible neither by $1+\tau z$, nor by $1-\tau z$ for any $\tau\ne0$, we reduce \eqref{e-oe5}
(possibly replacing $\tau$ by $-\tau$) to the conditions: for $k=1,2$,
$$F_k(z)=1+\sum_{i=1}^na_iz^{2i}+(-1)^k2z\quad
\text{is divisible by}\quad(1+(-1)^{k+1}\tau z)^{n+1},\quad \tau\ne0\ ,$$
which are equivalent to the following combinations:
\begin{equation}\begin{cases}&F_k(z),F_k'(z)\ \text{are divisible by}\ 1+(-1)^{k+1}\tau z,\\
&F_k''(z)\ \text{is divisible by}\ (1+(-1)^{k+1}\tau z)^{n-1},\end{cases}\quad\tau\ne0\ .\label{e-oe3}\end{equation}
The latter relations in \eqref{e-oe3} yield that
$$F_1''(z)=F''_2(z)=\sum_{i=1}^n2i(2i-1)a_iz^{2i-2}\quad \text{is divisible by}\quad (1-\tau^2z^2)^{n-1}\ ,$$
which results in
\begin{equation}a_i=\frac{(-1)^{i-1}}{i(2i-1)}\binom{i-1}{n-1}a_1\tau^{2i-2},\quad i=2,...,n\ ,
\label{e-oe4}\end{equation} whereas the former relations
in \eqref{e-oe3} lead to the system
\begin{equation} -2+\frac{a_1}{\tau}\sum_{i=1}^n\frac{2(-1)^i}{2i-1}\binom{i-1}{n-1}=1+\frac{2}{\tau}
+\frac{a_1}{\tau^2}\sum_{i=1}^n\frac{(-1)^{i-1}}{i(2i-1)}\binom{i-1}{n-1}=0\ .\label{e-oe6}\end{equation}
Observe that
$$\sum_{i=1}^n\frac{2(-1)^i}{2i-1}\binom{i-1}{n-1}=-2I_1,\quad I_1=\int_0^1(1-\xi^2)^{n-1}d\xi\ ,$$
where $I_1>0$, and
$$\sum_{i=1}^n\frac{(-1)^{i-1}}{i(2i-1)}\binom{i-1}{n-1}=2I_2,
\quad I_2=\int_0^1\left(\int_0^\xi(1-\eta^2)^{n-1}d\eta\right)d\xi\ ,$$
where $0<I_2<I_1$. Hence system \eqref{e-oe6} has a unique solution, which in view of
\eqref{e-oe4} implies the uniqueness of the sought curve $C$.

The last step in the proof of Theorem \ref{t-oe1} is to show that $C$ has no other singularity in $\C^2$, or, equivalently, that the polynomial $\Delta(x)$ has no multiple root other than $1/\tau^2$ if $(a_1,\dots, a_n)$ is the solution of equations \eqref{e-oe4} and \eqref{e-oe6}.
Indeed, let $x=\theta^2$, $\theta^2\neq 1/\tau^2$, be a root of $\Delta(x)$ of multiplicity $m\geq 2$. Then $z=(-1)^k\theta$ is a root of polynomial $F_k(z)$ of multiplicity $m$. Consider the curve $C'$ given by the equation $w-1-\sum_{i=1}^na_iz^{2i}=0$. The curve $C'$ is rational and it has singularity of type $T_{2n-1,2n}$ at infinity. In the pencil of lines $\lambda w+\mu z=0$, there are  at least three lines which intersect $C'$ in the number of points less than $2n=\deg C'$:
by assumption, each of the lines given by $w=2z$ and $w=-2z$ is tangent to $C'$ at least at two points with multiplicities $n+1$ and $m$, respectively, and the line given by $z=0$ intersects $C'$ at its singular point with multiplicity $2n-1$. Let $\overline C'$ be the normalization of $C'$. By Hurwits formula, applied to morphism $f:\overline C'\to\PP^1$ of degree $2n$ defined by the pencil of lines $\lambda w+\mu z=0$, we get the following inequality:
$$-2\geq -4n+ 2(n+m-1)+(2n-2)$$
which breaks down for $m\geq 2$.
\proofend

Consider the Fermat curve $C_n\subset \PP^2$ of degree $n\geq 3$ given by equation $z_1^n+z_2^n+z_3^n=0$.
By Pl\"ucker formulas, the dual curve $\hat C$ to $C$ has degree $n(n-1)$ and the following type of singularities:
$${\bf S}(\hat C_n)= 3nT_{n-1,n} +{\bf S}_{F_n},$$
where ${\bf S}_{F_n}$ is a sum of singularity types of the simple singularities of the curve dual to the Fermat curve $C_n$.

The following theorem will be used in the proof of Theorem \ref{rigit} (see subsection \ref{sec6}).
\begin{thm} \label{rigi} Let ${\bf d}_n=(n(n-1),1,1,1)$, ${\bf g}_n=(\frac{(n-1)(n-2)}{2},0,0,0)$,
and $n_{(1,2)}=n_{(1,3)}=n_{(1,4)}=n$,  $n_{(1)}=n_{(2,3)}=n_{(2,4)}=n_{(3,4)}=1$. Then the
family
$$V_n=V({\bf d}_n;{\bf g}_n;(n_{(1)}{\bf S}_{F_n}+(n_{(1,2)}+n_{(1,3)}+
n_{(1,4)})T^{n-1}_{n-1,n}+(n_{(2,3)}+n_{(2,4)}+n_{(3,4)})A_1)$$ is strictly rigid for each $n\geq 3$.
\end{thm}
{\bf Proof.} Consider a curve $\overline C=C_1\cup C_2\cup C_3\cup C_4\in V_n$ and denote the components $C_i$ of $\overline C$
as follows: $C_1=C$ and $C_i=L_{i-1}$ for $i\geq 2$. The curves $L_1$, $L_2$, $L_3$ are lines and $C$ is a curve of degree $n(n-1)$. The curve $C$ has $3n$ singular
points of the singularity type $T_{n-1,n}$, since $n_{(1,2)}=n_{(1,3)}=n_{(1,4)}=n$, and the set of all other singularities of $C$ is
${\bf S}_{F_n}$, that is, the singularity type of $C$ is the same as the singularity type of the dual curve to the Fermat curve of degree $n$.
Therefore, the dual curve $\hat C$ to $C$ is  non-singular and $\deg \hat C=n$, and to prove Theorem, it suffices to show that there are homogeneous coordinates in $\PP^2$ such that in these coordinate system the curve $\hat C$ is the Fermat curve.

Let $(z_1:z_2:z_3)$ be homogeneous coordinates in $\hat{\PP}^2$ such that
$p_1=(1:0:0)$, $p_2=(0:1:0)$, and $p_3=(0:0:1)$ are the points in $\hat{\PP}^2$ dual, respectively, to the
lines $L_1$, $L_2$, and $L_3$, and
let $$F(z_1,z_2,z_3)=\sum_{i_1+i_2+i_3=n}a_{\overline i}z^{\overline i}$$
be an equation of the curve $\hat C$, where $\overline i=(i_1,i_2,i_3)$ and $z^{\overline i}=z_1^{i_1}z_2^{i_2}z_3^{i_3}$.
Without loss of generality, we can assume that $a_{n,0,0}=a_{0,n,0}=a_{0,0,n}=1$ since the points $p_1$, $p_2$, and $p_3$
does not belong to the curve $\hat C$.

The curve $\hat C$ has $3n$ flex points given in local coordinates by equation $y=x^n$
and for $i=1,2,3$ there are $n$ lines $L_{i,1},\dots, L_{i,n}$ from the pencil of lines passing through $p_i$
such that they are tangent to $\hat C$ at its flex points.

Let $q_j=L_{1,j}\cap \hat C=(q_{1,j}:q_{2,j}:q_{3,j})$, $1\leq j\leq n$, and
let us rewrite the equation of the curve $\hat C$ in the form
\begin{equation} \label {eq-pr} F(z_1,z_2,z_3)=z_1^n+ \sum_{i=1}^{n}\binom{n}{i}
H_{i}(z_2,z_3)z_1^{n-i},
\end{equation}
where $H_i(z_2,z_3)$ are homogeneous polynomials in $z_2$, $z_3$ of degree $i$.
Then for $i=1,\dots,n-1$ the homogeneous polynomial $H_1^i(z_2,z_3)-H_i(z_2,z_3)$ of degree $i$ has
$n>i$ different roots, namely, $(q_{2,1}:q_{3,1}),\dots, (q_{2,n}:q_{3,n})$. Therefore
$H_i(z_2,z_3)=H_1^i(z_2,z_3)$ and hence the polynomial $F(z_1,z_2,z_3)$ has the form
\begin{equation} \label{eq-pr1}
F(z_1,z_2,z_3)=(z_1+H_{1}(z_2,z_3))^{n}+\widetilde H_n(z_2,z_3),
\end{equation}
where $\widetilde H_n(z_2,z_3)$ is a homogeneous polynomial of degree $n$,
$$\widetilde H_n(z_2,z_3)=\alpha\prod_{j=1}^n(q_{3,j}z_2-q_{2,j}z_3)$$
with some $\alpha\in\mathbb C$.

By the same arguments, we have
\begin{equation} \label{eq-pr2}
F(z_1,z_2,z_3)=(z_2+G_{1}(z_1,z_3))^{n}+\widetilde G_n(z_1,z_3)
\end{equation}
and
\begin{equation} \label{eq-pr10}
F(z_1,z_2,z_3)=(z_3+P_{1}(z_1,z_2))^{n}+\widetilde P_n(z_1,z_2),
\end{equation}
where $G_1(z_1,z_3)$, $P_1(z_1,z_2)$ and $\widetilde G_n(z_1,z_3)$, $\widetilde P_n(z_1,z_2)$ are homogeneous polynomials of
degree one and $n$ respectively.

It follows from \eqref{eq-pr2} that $\widetilde H_n(z_2,z_3)=(z_2+G_{1}(0,z_3))^{n}+\widetilde G_n(0,z_3)-H_1^n(z_2,z_3)$.
Therefore
\begin{equation} \label{eq-pr3}
F(z_1,z_2,z_3)=(z_1+az_2+bz_3)^n-(az_2+bz_3)^n+(z_2+cz_3)^n+(1-c^n)z_3^n
\end{equation}
for some $a,b,c\in\C$ (remind that, by assumption, $a_{n,0,0}=a_{0,n,0}=a_{0,0,n}=1$).

It follows from \eqref{eq-pr2} and \eqref{eq-pr3} that for $1\leq  k\leq n-1$
\begin{equation} \label{eq-pr4} G^k_1(z_1,z_3)=a^{n-k}((z_1+bz_3)^k-b^kz_3^k)+c^kz_3^k .\end{equation}
In particular,
$G_1(z_1,z_3)=a^{n-1}z_1+cz_3$  and hence
\begin{equation} \label{eq-pr5} (a^{n-1}z_1+cz_3)^k=a^{n-k}((z_1+bz_3)^k-b^kz_3^k)+c^kz_3^k .\end{equation}
It follows from \eqref{eq-pr5} that
\begin{equation} \label{eq-pr6} a^{j(n-1)}c^{k-j}=a^{n-k}b^{k-j} \end{equation}
for $1\leq  k\leq n-1$ and $1\leq j\leq k$. In particular, if we put $j=k$ in \eqref{eq-pr6}, then
we obtain that $a^n=1$ if $a\neq 0$.

Similarly, it follows from \eqref{eq-pr10} and \eqref{eq-pr3} that for $1\leq  k\leq n-1$
\begin{equation} \label{eq-pr7} P^k_1(z_1,z_2)=b^{n-k}((z_1+az_2)^k-a^kz_2^k)+c^{n-k}z_2^k .\end{equation}
In particular,
$P_1(z_1,z_3)=b^{n-1}z_1+cz_2$  and hence
\begin{equation} \label{eq-pr8} (b^{n-1}z_1+cz_2)^k=b^{n-k}((z_1+az_2)^k-a^kz_2^k)+c^{n-k}z_2^k .\end{equation}
It follows from \eqref{eq-pr8} that
\begin{equation} \label{last} \begin{array}{rll} b^{j(n-1)}c^{k-j} & = & b^{n-k}a^{k-j}, \\ c^k & =& c^{n-k}
\end{array} \end{equation}
for $1\leq  k\leq n-1$ and $1\leq j\leq k$. In particular, if $c\neq 0$ then we obtain that $c^2=c^n=1$, that is,
$c=\pm 1$ and $n$ is an even number if $c=-1$. If we put $j=k$ in \eqref{last}, then we obtain that $b^n=1$ if $b\neq 0$.

Let us show that the case when $abc\neq 0$ is impossible.
Indeed, if $abc\neq 0$, then $a^n=b^n=1$, $c=\pm 1$, and $c=-1$ only if $n$ is even. If we apply again \eqref{eq-pr6}
and \eqref{last} we obtain that
$c^{k-j}=(\frac{a}{b})^{k-j}=(\frac{b}{a})^{k-j}$ for $1\leq  k\leq n-1$ and $1\leq j\leq k$.
Therefore, $c=\frac{a}{b}=\frac{b}{a}$ and hence
$$F(z_1,z_2,z_3)=(z_1+az_2+acz_3)^n-(az_2+acz_3)^n+(z_2+cz_3)^n=(z_1+az_2+acz_3)^n,$$
since $a^n=c^n=1$. This contradicts the assumption that $C$ is an irreducible reduced curve.

It easily follows from \eqref{eq-pr6} and \eqref{last} that the case,  when the only one number
either $a$ or $b$, or $c$ is equal to zero, is also impossible.

In the case when $a=b=0$ and $c\neq 0$ we have
$F(z_1,z_2,z_3)=z_1^n+(z_2\pm z_3)^n$. But, it is impossible since $C$ is an irreducible curve.

The cases $a=c=0$, $b\neq 0$ and $b=c=0$, $a\neq 0$ are also impossible since in these cases we have,
respectively, that $F(z_1,z_2,z_3)=(z_1+bz_3)^n+ z_2^n$  since $b^n=1$ or $F(z_1,z_2,z_3)=(z_1+az_2)^n+ z_3^n$  since $a^n=1$.
As a result, we obtain that $a=b=c=0$, that is $F(z_1,z_2,z_3)=z_1^n+ z_2^n+z_3^n$. \proofend

\section{$k$-Rigid curves with $k\geq 2$}
\subsection{$2$-Rigid irreducible families of equisingular plane curves of degree $\geq 5$} \label {sec5}
An infinite series of irreducible $2$-rigid families is given in the following two theorems.
\begin{thm} \label{add2}
The families $V(2n+1;0;T_{n+1,2n+1}+T_{n,2n+1})$ are $2$-rigid and irreducible for $n\geq 2$.
\end{thm}
{\bf Proof.} Let $C\in V(2n+1;0;T_{n+1,2n+1}+T_{n,2n+1})$ and let $f:\PP^1\to\PP^2$ be a morphism such that $C=f(\PP^1)$. We can choose
non-homogeneous coordinate $t$ in $\PP^1$ and homogeneous coordinates $(z_1:z_2:z_3)$
in $\PP^2$ such that $f(t=0)=p_1=(0:0:1)$ is the singular point of $C$ of type $T_{n+1,2n+1}$, the line $z_1=0$ is the tangent line to $C$ at $p_1$,  $f(t=\infty)=p_2=(1:0:0)$ is the singular point of $C$ of type $T_{n,2n+1}$, the line $L_3$ given by $z_3=0$ is the
tangent line to $C$ at $p_2$. Then (maybe, after
change of the coordinates in $\PP^2$) the morphism $f$ is given by
\begin{equation} \label{par} z_1=t^{2n+1}, \, \,
z_2=t^{n+1}, \, \, z_3=at+1\end{equation}
for some $a\in \C$.
Let $C_a$ has parametrization \eqref{par}.
The intersection number of $C_a$ and $L_3$ at the point $p_2$ is
$$(C_0,L_3)_{p_2}=\left\{ \begin{array}{ll} 2n+1\,\, & \text{if}\, a=0\\
2n\,\, & \text{if}\, a\neq 0. \end{array}\right.  $$
Hence there is no projective transformation $h\in \Aut(\PP^2)$ such that $h(C_0)=C_1$.
On the other hand, if $a\neq 0$, we make change
$t_1=at$, and then the projective transformation
$h((z_1:z_2:z_3))=(a^{2n+1}z_1:a^{n+1}z_2:z_3)$
sends  $C_1$ to $C_a$. \proofend

\begin{thm} \label{add3}
The families $V(4n;0;T_{2n-1,4n}+T_{2n+1,4n})$ are $2$-rigid and irreducible for $n\geq 3$.
\end{thm}
{\bf Proof.} Let $C\in V(4n;0;T_{2n-1,4n}+T_{2n+1,4n})$ and let $f:\PP^1\to\PP^2$ be a
morphism such that $C=f(\PP^1)$. As in the proof of Theorem \ref{add2}, it is easy to show that  we can choose
non-homogeneous coordinate $t$ in $\PP^1$ and homogeneous coordinates $(z_1:z_2:z_3)$
in $\PP^2$ such that $f(t=0)=p_1=(0:0:1)$ is the singular point of $C$ of type $T_{2n-1,4n}$,
the line $z_1=0$ is the tangent line to $C$ at $p_1$,  $f(t=\infty)=p_2=(1:0:0)$ is the singular
point of $C$ of type $T_{2n+1,4n}$, and the line $z_3=0$ is the
tangent line to $C$ at $p_2$. Then (maybe, after
change of the coordinates in $\PP^2$ and $\PP^1$) the morphism $f$ is given by
\begin{equation}\label{par1} z_1=(t^2-a^2)t^{4n-2}, \, \,
z_2=t^{2n-1}, \, \,z_3=1 \end{equation}
for some $a\in \C$.
Let $C_a$ has parametrization \eqref{par1}.
The intersection number of $C_a$ and $L_1$ at the point $p_1$ is
$$(C_0,L_1)_{p_1}=\left\{ \begin{array}{ll} 4n\,\, & \text{if}\, a=0\\
4n-2\,\, & \text{if}\, a\neq 0. \end{array}\right.  $$
Hence there is no projective transformation $h\in \Aut(\PP^2)$ such that $h(C_0)=C_1$.
On the other hand, if $a\neq 0$, we make change
$t_1=a^{-1}t$, and then the projective transformation
$h((z_1:z_2:z_3))=(a^{-4n}z_1:a^{1-2n}z_2:z_3)$
sends  $C_a$ to $C_1$.
\proofend

\subsection{$k$-Rigid families parameterizing curves with $k$ connected components} \label {sec6}
\begin{thm} \label{rigit} For each $k\in \N$ there is a $k$-rigid family of equisingular plane curves consisting of $k$ irreducible components.
\end{thm}
{\bf Proof.} Let $n=2k+1$ and let $n_{(1,2)}=n-2$, $n_{(1,3)}=n_{(1,4)}=n$, and $n_{(1)}=n_{(2,4)}=n_{(3,4)}=n_{(1,2,5)}=n_{(1,2,6)}=n_{(3,4,5,6)}=1$.
Consider the family
$\overline V_n=V(\overline {\bf d}_n;\overline {\bf g}_n;\overline {\bf S}_n)$, where
 $$\overline {\bf d}_n=(n(n-1),1,1,1,1,1),\qquad \overline {\bf g}_n=(\frac{(n-1)(n-2)}{2},0,0,0,0,0),$$ and
$$\begin{array}{ll} \overline {\bf S}_n = & (n_{(1,2)}+n_{(1,3)}+
n_{(1,4)})T^{n-1}_{n-1,n}+(n_{(1,2,5)}+n_{(1,2,6)})T^n_{n-1,n} \\ & +n_{(1)}{\bf S}_{F_n}+n_{(3,4,5,6)}T_{4,4}+(n_{(2,4)}+n_{(3,4)})A_1),\end{array}$$
where ${\bf S}_{F_n}$ is a sum of singularity types of the simple singularities of the curve dual to the Fermat curve given by equation $z_1^n+z_2^n+z_3^n=0$.

Consider a curve $\widetilde C=\bigcup_{i=1}^6C_i\in \overline V_n$.
Denote the curve $C_1$ by $C$ and the curve $C_i$ by $L_{i-1}$ for $i\geq 2$. The curves $L_1, \dots, L_5$
are lines and $C$ is a curve of degree $n(n-1)$ and it is easy to see that
$\overline C=C\cup L_1\cup L_2\cup L_3\in V_n$, where $V_n$ is the
family of plane curves from Theorem \ref{rigi}.
By Theorem \ref{rigi}, we can assume that $C$ is the curve dual to the Fermat curve of degree $n$ and
$L_i$ is given by equation $z_i=0$ for $i=1,2,3$. Then it follows from the singularity type of the curve
$\widetilde C$ that $L_4$ and $L_5$ have, respectively, equations $z_2+\varepsilon^{m_1}z_3=0$ and
$z_2+\varepsilon^{m_2}z_3=0$, where $m_1\not\equiv m_2(mod\, n)$ and $\varepsilon$ is a primitive root of
the equation $x^n+1=0$.

Consider two curves $\widetilde C_i=\overline C\cup L_{4,i}\cup L_{5,i}\in \overline V_n$, $i=1,2$,
where $L_{4,i}$ is given by equation $z_2+\varepsilon^{m_{1,i}}z_3=0$ and $L_{5,i}$ is given by equation $z_2+\varepsilon^{m_{2,i}}z_3=0$.
It is easy to see that a projective transformation $h\in \Aut(\PP^2)$ such that $\widetilde C_2=h(\widetilde C_1)$ exists
if and only if $m_{1,1}-m_{1,2}\equiv \pm(m_{2,1}-m_{22})
\mod n$. Hence
$\overline V_n$ is a $(\frac{n-1}{2})$-rigid family of plane curves consisting of $\frac{n-1}{2}=k$ irreducible components.  \proofend

{\ncsc Steklov Mathematical Institute \\[-21pt]

Gubkina str., 8 \\[-21pt]

Moscow, Russia. \\[-21pt]

{\it E-mail address}: {\ntt kulikov@mi.ras.ru}

\vskip10pt

{\ncsc School of Mathematical Sciences \\[-21pt]

Raymond and Beverly Sackler Faculty of Exact Sciences\\[-21pt]

Tel Aviv University \\[-21pt]

Ramat Aviv, 69978 Tel Aviv, Israel} \\[-21pt]

{\it E-mail address}: {\ntt shustin@post.tau.ac.il}


\begin{thebibliography}{99}

\bibitem{Be} G. V. Belyi. On Galois extensions of a maximal cyclotomic field
{\it Izv. Akad. Nauk SSSR, Ser. Mat.} {\bf 43} (1979), no. 2, 267--276.

\bibitem{GLS} G.-M. Greuel, C. Lossen, and E. Shustin.
Equisingular families of projective curves. In: {\it Global aspects of complex geometry}/F. Catanese et al., eds., Springer, 2006, pp. 171--209.

\bibitem{Hir} F. Hirzebruch. Some examples of algebraic surfaces. {\it Contemp. Math.} {\bf 9} (1982), 55--71.

\bibitem{Ku} Vik. S. Kulikov. {\it A Remark on Classical Pluecker's formulae}.
Preprint at arXiv:1101.5042, (2011).

\bibitem{K-H} Vik. S. Kulikov, V. M. Kharlamov. On real structures on rigid surfaces.
{\it Izv. Math.} {\bf 66:1} (2002), 133-–150.

\bibitem{Pa} K. H. Paranjape. A geometric characterization of arithmetic varieties.
{\it Proc. Indian Acad. Sci. Math. Sci.} {\bf 112} (2002),  no. 3, 383--391.

\bibitem{Sh} E. I. Shustin. Versal deformations in the space of plane curves of fixed degree.
{\it Func. Anal. Appl.} {\bf 21} (1987), 82--84.

\bibitem{ZO} M. G. Zaidenberg and S. Yu. Orevkov.
On rigid rational cuspidal plane curves.
{\it Russ. Math. Surveys} {\bf 51} (1996), no. 1, 179--180.

\bibitem{Z1} O. Zariski. Studies in singularity. I. Equivalent singularities of plane algebroid curves,
{\it Amer. J. Math.} {\bf 87} (1965), 507--536.

\bibitem{Z2} O. Zariski. Studies in singularity. II. Equisingularity in codimension 1 (and characteristic zero),
{\it Amer. J. Math.} {\bf 87} (1965), 972--1006.

\bibitem{Wal} C.T.C. Wall. Geometry of quartic curves. {\it Math. Proc. Camb. Phil. Soc.} {\bf 117} (1995), 415--423.

\end{thebibliography}
\end{document}